\documentclass[12pt]{article}
\usepackage{amsmath}
\usepackage{amssymb}

\usepackage[cp1251]{inputenc}
\usepackage[russian]{babel}
\newcommand{\il}[2]{\int\limits_{#1}^{#2}}

\newcommand{\ph}{\phantom{a}}
\newcommand{\phh}{\phantom{aaa}}

\newcommand{\sist}[2]{\left\{
\begin{array}{l}
{#1}\\
\ph\\
{#2}
\end{array}
\right.}

\begin{document}
MSC  34L30,   34C99

\vskip 10pt
{\bf \centerline {Periodic solutions of  quaternionic}

\centerline{Riccati equations with periodic coefficients}

\vskip 20pt \centerline {G. A. Grigorian}}
\centerline{\it Institute  of Mathematics of NAS of Armenia}
\centerline{\it E -mail: mathphys2@instmath.sci.am}
\vskip 20 pt

Abstract. In this paper we study the conditions, under which the quaternionic Riccati equations have periodic solutions. The obtained result we compare with one recently obtained important one.


\vskip 20pt
{\bf 1. Introduction}. Let  $a(t),\phantom{a} b(t),\phantom{a} c(t)$ and $d(t)$  be quaternionic-valued  continuous functions on $[t_0,+\infty)$, i.e.,
$a(t)\equiv a_0(t) + i a_1(t) + j a_2(t) + k a_3(t), \phantom{a} b(t)\equiv b_0(t) + i b_1(t) + j b_2(t) + k b_3(t), \phantom{a} c(t)\equiv c_0(t) + i c_1(t) + j c_2(t) + k c_3(t), \phantom{a} d(t)\equiv d_0(t) + i d_1(t) + j d_2(t) + k d_3(t),$
where  $a_n(t), \phantom{a} b_n(t), \phantom{a} c_n(t), \phantom{a} d_n(t) \phantom{a} (n=\overline{0,3})$ are
 real-valued continuous functions on  $[0,+\infty)$, $i,\phantom{a} j, \phantom{a} k$ are the imaginary unities, satisfying the conditions
$$
i^2 = j^2 = k^2 = ijk = -1, \phantom{a} ij = - ji = k.
$$
Consider the  Riccati equation
$$
q' + q a(t) q + b(t) q + q c(t) + d(t) = 0, \phantom{aaa} t\ge 0. \eqno (1.1)
$$
This equation appear in various problems of mathematics, in particular, in problems of mathematical physics (e. g., in the Euler's vorticity dynamics [1], in the Euler's fluid dynamics [2], in the problem of classification of diffeomorphisms of $\mathbb{S}^4$ [3], and in  other areas of the natural science (see, e.g.  [4, 5] and cited works therein). In the past years the problem of studying the question of finding the conditions under which quaterinionic Riccati equations with periodic coefficients have periodic solutions attract the attention of many authors (see [1-5], and cited works therein). In the work  [2] it was shown that  for enough small norm $\max\limits_{t\in[0,T]}|d(t)|$ the equation
$$
q' + q^2 + d(t) = 0  \eqno (1.2)
$$
with the $T$-periodic $d(t)$ has at least one $T$-periodic solution. In the mentioned work it was shown also that there exists a real-valued $d(t)$ for which Eq. (1.2) has no $T$-periodic solutions.
Among the results from this direction notice the following important result of S. P. Wilzinski.

{\bf Theorem 1.1 [1, Theorem 3 and Remark 4]}. {\it Let the coefficients $a, \ph b, \ph c, \ph d \in \mathbb{C}(\mathbb{R}, \mathbb{H})$ be T-periodic and $Ark(0) = 0$. If the conditions

\noindent
(i) \phh $a d \not \equiv 0, \phh Ark [a] < \frac{\pi}{4},$

\noindent
(ii) \phh $\max\limits_{t \ge 0} Ark[a(t)] + \max\limits_{t \ge 0} Ark[-d(t)] \le \frac {\pi}{2}$,

\noindent
(iii) \phh $Re \hskip 2pt [b + c] \le 0, \phh Im \hskip 2pt [b + c] \equiv 0$

\noindent
hold, then Eq. (1.1) ha exactly two periodic solutions $\xi, \eta$ in $\mathbb{H}$. Moreover $Re \hskip 2pt [\eta] < 0$ and $\eta$ is asymptotically stable, while $Re \hskip 2pt [\xi] < 0$ and $\xi$ is asymptotically unstable, and every other solution in $\mathbb{S}^4$ is holomorphic  to them. Every non periodic solution, starting $\mathbb{S}(-\frac{\pi}{2})$ is f. b. or enters $-\mathbb{S}(-\frac{\pi}{2})$. Every solution, starting in $-\mathbb{S}(-\frac{\pi}{2})$ slays on for all $t \ge 0$. There are no b. f. b. solution.}

Here $\mathbb{H}$ denotes the algebra of quaternions, $Ark [q] \equiv |Arg (s_q + i|v_q||, \ph s_q \equiv q_0 = Re \hskip 2pt q, \ph v_q\equiv iq_1 + j q_2 + k q_3 = Im \hskip 2pt q$ for any quaternion $q\equiv q_0 + i q_1 + j q_2 + k q_3$, f. b. forward blowing up, b. b. backward blowing up,  b. f. b. backward and forward blowing up. This result is obtained by the use of topological and geometrical methods. Unlike of these methods in this paper we use only analitical methods for obtain a criterion of existence of a periodic solution for Eq. (1.1) with periodic coefficients. The obtained result we compare with Theorem 1.1.

\hskip 10 pt

{\bf 2. Auxiliary propositions}. Set: $p_{0,m}(t)\equiv b_m(t) + c_m(t), \ph m=\overline{1,3}, \ph p_{1,1}(t) \equiv b_1(t) + c_1(t), \ph p_{1,2}(t) \equiv b_2(t) - c_2(t), \ph p_{1,3}(t) \equiv b_3(t) - c_3(t), \ph p_{2,1}(t) \equiv b_1(t) - c_1(t), \ph p_{2,2}(t) \equiv b_2(t) + c_2(t), \ph p_{2,3}(t) \equiv b_3(t) - c_3(t), \ph p_{3,m}(t) \equiv b_m(t) - c_m(t(), \ph m=\overline{1,3},$
$$
D_0(t)\equiv \sist{\sum\limits_{m=1}{3}p_{0,m}^2(t) + 4 a_0(t)d_0(t), \ph if \ph a_0(t) \ne 0,}{4 d_0(t), \ph if \ph a_0(t) = 0,} \phantom{aaaaaaaaaaaaaaaaaaaaaaaaaaaaa}
$$
$$
\phantom{aaaaaaaaaaaaaaaaaaaaaaaaaaaaaa}D_n(t)\equiv \sist{\sum\limits_{m=1}{3}p_{n,m}^2(t) - 4 a_0(t)d_0(t), \ph if \ph a_0(t) \ne 0,}{-4 d_n(t), \ph if \ph a_0(t) = 0,}
$$
$n = \overline{1.3}, \ph t \ge 0$.
Let $\mathfrak{S}$ be a nonempty  subset of the set $\{0, 1, 2, 3\}$ and let $\mathfrak{O}$ be its complement, i.e., $\mathfrak{O} = \{0, 1, 2, 3\} \backslash \mathfrak{S}$

{\bf Theorem 2.1 [6, Theorem 3.1]}. {\it Assume  $a_n(t) \ge 0, \phantom{a}  n \in \mathfrak{S}$ and if $a_n(t) = 0$  then $p_{n,m}(t) = 0, \phantom{a} m =~ \overline{1,3}, \ph   n\in \mathfrak{S}; \phantom{a} a_n(t) \equiv 0, \phantom{a}  n \in  \mathfrak{O}, \phantom{a} D_n(t) \le 0, \phantom{a} n \in \mathfrak{S}, \phantom{a} t \ge t_0.$
Then for every  $\gamma_n \ge 0, \phantom{a} n \in \mathfrak{S}, \phantom{a} \gamma_n \in (- \infty; + \infty), \phantom{a} n \in \mathfrak{O}$,  Eq. (1.2) has a solution  $q(t) \equiv q_0(t) - i q_1(t) - j q_2(t) - k q_3(t)$ on $[t_0; +\infty)$  with  $q_n(t_0) = \gamma_n, \phantom{a} n = \overline{0,3}$  and
$$
q_n(t) \ge 0, \phantom{aaa} n \in \mathfrak{S}, \phantom{aaa}   t\ge t_0.
$$
Moreover if for some $n \in \mathfrak{S}$  $\gamma_n > 0$,  then also  $q_n(t) > 0$.}

{\bf Definition 2.1.} {\it A solution $q(t)$ of Eq. (1.1), existing on $[0,+\infty)$,  is called $t_1$-normal ($t_1 \ge 0$) if there exists a neighborhood $U(q(t_1))$ of $q(t_1)$ such that every solution $\widetilde{q}(t)$ of Eq. (1.1) with $\widetilde{q}(t_1) \in U(q(t_1))$ exists on $[t_1,+\infty)$, otherwise $q(t)$ is called $t_1$-extremal.}

For any solutions $q_1(t)$ and $q_2(t)$ of Eq. (1.1) existing on $[t_1, +\infty) \ph (t_1 \ge 0)$ set
$$
\mathbb{I}_{q_1,q_2}(t) \equiv \il{t_1}{t} Re \hskip 2pt[a(\tau) (q_1(\tau) - q_2(\tau))]d \tau, \ph t \ge t_1.
$$

{\bf Theorem 2.2 [7, Corollary 3.1]}. {\it The following statements are valid:

\noindent
1) any two solutions   $q_1(t)$ and $q_2(t)$ of Eq. (1.1), existing on $[t_1,+\infty)$,  are $t_1$-normal if and only if the function $\mathbb{I}_{q_1,q_2}(t)$
is bounded;

\noindent
2) if $q_N(t)$ and $q_*(t)$ are $t_1$-normal and $t_1$-extremal solutions of Eq. (1.2) respectively then
$$
\limsup\limits_{t \to + \infty} \il{t_1}{t} Re [a(\tau)(q_*(\tau) - q_N(\tau))] d \tau < +\infty,
$$
$$
\liminf\limits_{t \to + \infty} \il{t_1}{t} Re [a(\tau)(q_*(\tau) - q_N(\tau))] d \tau = -\infty;
$$

$
\phantom{aaaaaaaaaaaaaaaaaaaaaaaaaaaaaaaaaaaaaaaaaaaaaaaaaaaaaaaaa} \blacksquare
$}

\hskip 10 pt

{\bf 3. Main result.} For any continuous functions on $[0,+\infty)$ set
$$
\biggl[\frac{u(t)}{v(t)}\biggr]_0\equiv\sist{\frac{u(t)}{v(t)}, \ph if \ph v(t)\ne 0,}{0, \phh if \ph v(t) = 0, \ph t \ge 0.}
$$

{\bf Theorem 3.1.} {\it Let $a(t), \ph b(t), \ph c(t)$ and $d(t)$ be $T$-periodic and let the following \linebreak conditions be satisfied.

\noindent
1) $a_n(t) \ge 0, \ph n=0,1, \ph a_0(t) + a_1(t) \not \equiv 0, \ph t \ge 0, \ph a_n(t)\equiv 0, \ph n=2,3,$

\noindent
2) $ D_n(t)\le 0, \ph t \ge 0, \ph D_n(t) \not \equiv 0, \ph n=1,2$,

\noindent
3) $supp \hskip 2pt [b_n(t) + c_n(t)] \subset supp \hskip 2pt [a_0(t)], \ph supp \hskip 2pt [b_n(t) - c_n(t)] \subset   supp \hskip 2pt [a_1(t)]$,

\noindent
4) $\bigl[\frac{b_n(t) + c_n(t)}{a_0(t)}\bigr]_0$ and  $\bigl[\frac{b_n(t) - c_n(t)}{a_1(t)}\bigr]_0, \ph n =2,3$ are bounded on $[0,+\infty)$,

\noindent
5) $\il{0}{T}[b_0(\tau) + c_0(\tau)] d\tau \ge 0.$

Then there exists a natural number $m_0$ such that Eq. (1.1) has a $m_0 T$-periodic \linebreak solution $q(t)\equiv q_0(t) - i q_1(t) - j q_2(t) - k q_3(t)$ such that $q_n(t) \ge 0, \ph n = 0, 1, \ph t \ge 0$.
If in addition $D_0(t)D_1(t)\not\equiv 0$, then every solution $\widetilde{q}(t)\equiv \widetilde{q}_0(t) - i \widetilde{q}_1(t) - j \widetilde{q}_2(t) - k \widetilde{q}_3(t)$ of Eq. (1.1) with $\widetilde{q}_n(0) \ge 0, \ph n=0,1$ (in particular $q(t)$) is $0$-normal and the integral
$$
I_{q,\widetilde{q}}(t)\equiv \il{0}{t}Re \hskip 2pt [a(\tau)(q(\tau) - \widetilde{q}(\tau))]d \tau, \phh t \ge 0
$$
is bounded.
}

Proof. It follows from the condition 1) that Eq. (1.1) is equivalent to the system
$$
\left\{\begin{array}{l} q_0'+a_0(t) q_0^2+\{b_0(t) + c_0(t) + 2 a_1(t) q_1\} q_1 - a_0(t) q_1^2 + (b_1(t) + c_1(t) q_1 -\\  \phantom{aaaaaaaaaaaaaaaaaaaaaaaaaaaaaaaaaaaaaaaaaaaaaa} -\widetilde{P}(t,q_2,q_3) = 0,\\
q_1' + a_1(t) q_1^2 + \{b_0(t) + c_0(t) + 2 a_0(t) q_0\} q_0 - a_1(t) q_0^2 + (b_1(t) + c_1(t) q_0 -\\ \phantom{aaaaaaaaaaaaaaaaaaaaaaaaaaaaaaaaaaaaaaaaaaaaaa} -\widetilde{Q}(t,q_2,q_3) = 0,\\
q_2' + \{b_0(t) + c_0(t) + 2[a_0(t)q_0 + a_1(t) q_1]\} q_2 - (b_1(t) - c_(t)) q_3 \\ -\phantom{aaaaaaaaaaaaaaaaaaaaaaaaaaaaaaaaaaaaaaaaaaaa}  -\widetilde{R}(t,q_0,q_1)=0, \\
q_3' + \{b_0(t) + c_0(t) + 2[a_0(t)q_0 + a_1(t) q_1]\} q_3 - (b_1(t) + c_(t)) q_2 -\\ \phantom{aaaaaaaaaaaaaaaaaaaaaaaaaaaaaaaaaaaaaaaa} -\widetilde{S}(t,q_0,q_1)=0, \ph t \ge 0,
\end{array}
\right. \eqno (3.1)
$$
where $\widetilde{P}(t,q_2,q_3) \equiv a_0(t)[q_2^2 + q_3^2] +[b_2(t) + c_2(t)]q_2 + [b_3(t) + c_3(t)] q_3 + d_0(t),$

\phantom{aaaaa}$\widetilde{Q}(t,q_2,q_3) \equiv a_1(t)[q_2^2 + q_3^2] +[b_3(t) + c_3(t)]q_2 + [b_2(t) - c_2(t)] q_3 - d_1(t),$

\phantom{aaaaaaaa}$\widetilde{R}(t,q_0,q_1)\equiv (b_2(t) + c_2(t)) q_0  - (b_3(t) - c_3(t)) q_1 - d_2(t),$

\phantom{aaaaaaaaaaaa}$\widetilde{S}(t,q_0,q_1)\equiv (b_3(t) + c_3(t)) q_0  - (b_2(t) - c_2(t)) q_1 - d_3(t), \ph t \ge 0.$

Let $q(t)\equiv q_0(t) - i q_1(t) - j q_2(t) - k q_3(t) \ph (q_n(t) \in \mathbb{R}, \ph n = \overline{0,3}, \ph t \ge 0)$ be a solution of Eq. (1.1) with $q_n(0) \ge 0, \ph n= 0,1.$ By virtue of Theorem 2.1 from 1) and 2) it follows that $q(t)$ exists on $[0,+\infty)$ and
$$
q_n(t) \ge 0, \ph n=0,1, \ph t \ge 0. \eqno (3.2)
$$
Show that there exists a natural number $m_0$, not depending on $q(t)$, such that
$$
\il{0}{m_0T}[b_0(\tau) + c_0(\tau) + a_0(\tau) q_0(\tau) + a_1(\tau) q_1(\tau)] d \tau \ge 1. \eqno (3.3)
$$
From the condition 5) it follows that $I_0 \equiv \il{0}{T}[b_0(\tau) + c_0(\tau)] d \tau > 0$ or $I_0 = 0$. If $I_0 > 0$, then from 1) and (3.2) it follows that (3.3) is satisfied for $m_0 = \bigl[\frac{1}{I_0}\bigr] + 1$. Assume
$$
I_0 = 0. \eqno (3.4)
$$
By the first equation of the system (3.1) we can interpret $q_0(t)$ as a solution of the linear equation
$$
v'+\{b_0(t) + c_0(t) + a_0(t) q_0(t) + 2 a_1(t) q_1(t)\} v - a_0(t) q_1^2(t) + (b_1(t) + c_1(t) q_1(t) - \phantom{aaaaaaaaaaaaaaaa}
$$
$$
\phantom{aaaaaaaaaaaaaaaaaaaaaaaaaaaaaaaaaaaaaaaaaaaaaaaa}-\widetilde{P}(t,q_2(t),q_3(t)) = 0, \ph t \ge 0.
$$
Then according to the Cauchy formula we have
$$
q_0(t) = a_0(0)\exp\biggl\{-\il{0}{t}[a_0(\tau) q_0(\tau) + b_0(\tau) + c_0(\tau) + 2 a_1(\tau) q_1(\tau)]d \tau\biggr\} + \phantom{aaaaaaaaaaaaaaaa}
$$
$$
+\il{0}{t}\exp\biggl\{-\il{\tau}{t}[a_0(s) q_0(s) + b_0(s) + c_0(s) + 2 a_1(s) q_1(s)] d s\biggr\}\biggl[a_0(\tau) q_1^2(\tau)- (b_1(\tau) +
$$
$$
\phantom{aaaaaaaaaaaaaaaaaaaaaaaaaaaaaaa}+ c_1(\tau)) q_1(\tau) + \widetilde{P}(\tau, q_2(\tau),q_3(\tau))\biggr]d\tau, \phh t \ge 0.
$$
Since $a_0(t)q-1^2(t) - [b_1(t) + c_1(t)] q_1(t) + \widetilde{P}(t, q_2(t),q_3(t)) \ge - D_0(t), \ph t \ge 0$, from the last equality we obtain
$$
\zeta_{0,m}\equiv \il{0}{m T} a_0(\tau) q_0(\tau) d \tau \ge
$$
$$
 \ge \il{0}{m T}\Biggl(\il{0}{t}\exp\biggl\{-\il{\tau}{t}[a_0(s) q_0(s) + b_0(s) + c_0(s) + 2 a_1(s) q_1(s)] d s\biggr\}\bigr[-D_0(\tau)\bigr]d\tau\Biggr) d t,  \eqno (3.5)
$$
$m =1,2,\dots.$ Analogously on the basis of the second equation of the system (3.1) it can be obtained the estimates
$$
\zeta_{1,m}\equiv \il{0}{m T} a_1(\tau) q_1(\tau) d \tau \ge
$$
$$
 \ge \il{0}{m T}\Biggl(\il{0}{t}\exp\biggl\{-\il{\tau}{t}[a_1(s) q_1(s) + b_0(s) + c_0(s) + 2 a_0(s) q_0(s)] d s\biggr\}\bigr[-D_1(\tau)\bigr]d\tau\Biggr) d t,  \eqno (3.6)
$$
$m =1,2,\dots.$ Since by 1) and (3.2) $a_n(t) q_n(t) \ge 0, \ph n=0,1, \ph t \ge 0,$ we have
$$
\il{\tau}{t}a_n(s) q_n(s) d s \le \xi_{n,m}, \ph 0 \le \tau \le t \le m T,\ph n=0,1, \ph m=1,2,\dots. \eqno (3.7)
$$
It follows from (3.4) that $\exp\biggl\{\il{\tau}{t}(b_0(\tau) + c_0(\tau)) d \tau\biggr\} \ge \varepsilon_0, \ph 0 \le \tau \le t \le m T$ for some $\varepsilon > 0$ and for all $m = 1,2,\dots$.  This together with (3.5)-(3.7) implies that
$$
\xi_{0,m} \xi_{1,m} e^{3[\xi_{0,m} + \xi_{1,m}]} \ge \varepsilon_0^2\il{0}{mT}d t\il{0}{t}[-D_0(\tau)]d\tau \il{0}{mT}d t\il{0}{t}[-D_1(\tau)]d\tau, \ph m=1,2,\dots. \eqno (3.8)
$$
It follows from the condition 2) that the functions $f_n(t)\equiv \il{0}{t}[-D_n(\tau)]d\tau, \ph n=0,1$ are  non decreasing and $>0$ for all enough large $t > 0$. Then we can chose $m_0$ so large that
$$
\varepsilon_0^2\il{0}{mT}d t\il{0}{t}[-D_0(\tau)]d\tau \il{0}{mT}d t\il{0}{t}[-D_1(\tau)]d\tau \ge e^3.
$$
From here and from (3.8) it follows that
$$
\xi_{0,m} \xi_{1,m} e^{3[\xi_{0,m} + \xi_{1,m}]} \ge e^3. \eqno (3.9)
$$
From here it follows
$$
\xi_{0,m} + \xi_{1,m} \ge 1. \eqno (3.10)
$$
Indeed, if $\xi_{0,m} + \xi_{1,m} <1$, then since $\xi_{n,m_0} \ge 0, \ph n =0,1$ we have $\xi_{0,m_0}  \xi_{1,m_0} < 1$ and, therefore $\xi_{0,m} \xi_{1,m} e^{3[\xi_{0,m} + \xi_{1,m}]} < e^3$, which contradicts (3.9). From (3.9) and from nonnegativity of $a_n(t) q_n(t), \ph n=0,1$ on $[0,+\infty)$ it follows (3.3). Set: $z(t)\equiv q_0(t) -~ iq_1(t), \linebreak \widetilde{b}(t)\equiv b_0(t) + i b_1(t), \ph \widetilde{c}(t)\equiv c_0(t) + i c_1(t), \ph U(t)\equiv \widetilde{P}(t,q_2(t),q_3(t)) - i\widetilde{Q}(t,q_2(t),q_3(t)), \linebreak w(t)\equiv q_2(t) + i q_3(t), \ph A(t)\equiv b_0(t) + c_0(t) + 2[a_0(t) q_0(t) + a_1(t) q_1(t)] + i [b_1(t) - c_1(t)], \linebreak V(t) \equiv-(b_2(t) + c_2(t)) q_0(t) + (b_3(t) - c_3(t)) q_1(t) - i [(b_3(t) + c_3(t)) q_0(t) + (b_2(t) - c_2(t))] -d_2(t) - i d_3(t), \ph t \ge 0.$ Then by Cauchy formula and by (3.1) we obtain
$$
z(t) = \exp\biggl\{-\il{0}{t}[a_0(\tau) z(\tau) + \widetilde{b}(\tau) + \widetilde{c}(\tau)]d \tau\biggr\} z(0) + \phantom{aaaaaaaaaaaaaaaaaaaaaaaaaaaaaaaaaaaaaa}
$$
$$
\phantom{aaaaaaaaaaaaa}+\il{0}{t}\exp\biggl\{-\il{\tau}{t}[a(s) z(s) + \widetilde{b}(s) + \widetilde{c}(s)] d s\biggr\} U(\tau) d \tau, \ph t \ge 0. \eqno (3.11)
$$
$$
w(t) = \exp\biggl\{-\il{0}{t} A(\tau) d\tau\biggr\} w(0) - \il{0}{t}\exp\biggl\{-\il{\tau}{t}A(s) d s\biggr\} V(\tau) d\tau, \phh t \ge 0. \eqno (3.12)
$$
Consider the functions
$$
J(t)\equiv \il{0}{t}\exp\biggl\{-\il{\tau}{t}A(s) d s\biggr\} V(\tau) d\tau, \ph \widetilde{V}(t) \equiv V(t) + d_2(t) + i d_3(t) = - (b_2(t) + c_2(t)) q_0(t) +
$$
$$
+ (b_3(t) - c_3(t)) q_1(t) - i [(b_3(t) + c_3(t)) q_0(t) + (b_2(t) - c_2(t))q_1(t)], \ph t \ge 0.
$$
It follows from the conditions 3) and 4) that $\biggl[\frac{\widetilde{V}(t)}{2(a_0(t) q_0(t) + a_1(t) q_1(t))}\biggr]_0$ exists on $[0,m_0]$, is bounded on it and
$$
J(t) = \exp\biggl\{-2\il{0}{t}[a_0(\tau) q_0(\tau) + a_1(\tau) q_1(\tau)] d \tau\biggr\}\il{0}{t}\biggl[\exp\biggl\{2\il{0}{\tau}(a_0(s) q_0(s) + a_1(s) q_1(s)) d s\biggr\}\biggr]'\times
$$
$$
\times\exp\biggl\{-\il{\tau}{t}[b_0(s) + c_0(s) + i (b_1(s) - c_1(s))] d s\biggr\}\biggl[\frac{\widetilde{V}(\tau)}{2[a_0(\tau) q_0(\tau) + a_1(\tau) q_1(\tau)]}\biggr]_0 d \tau +
$$
$$
\phantom{aaaaaaaaaaaaa} +\il{0}{t}\exp\biggl\{-\il{\tau}{t} A(s) d s\biggr\}\bigl[-d_2(\tau) + id_3(\tau)\bigr] d \tau, \ph t \in [0,m_0 T]. \eqno (3.13)
$$
Obviously
$$
\biggl|\exp\biggl\{-\il{\tau}{t}\bigl[b_0(s) + c_0(s) + i(b_1(s) -c_1(s))\bigr]d s\biggr\}\biggl[\frac{\widetilde{V}(\tau)}{2(a_0(\tau) q_0(\tau) + a_(\tau) q_1(\tau))}\biggr]_0 \biggr| \le
$$
$$
\le \exp\biggl\{-\il{\tau}{t}(b_0(s) + c_0(s))ds\biggr\}\biggl[\frac{|b_2(\tau) + c_2(\tau)| + |b_3(\tau) + c_3(\tau)|}{2a_0(\tau)} +
$$
$$
+\frac{|b_2(\tau) - c_2(\tau)| + |b_3(\tau) - c_3(\tau)|}{2a_1(\tau)}\biggr]_0, \ph 0 \le \tau \le t \le m_0 T \eqno (3.14)
$$
and
$$
\biggl|\il{0}{t}\exp\biggl\{-\il{\tau}{t}A(s)d s\biggr\}\bigl[d_2(\tau) + i d_3(\tau)\bigr]d \tau\biggr| \le \phantom{aaaaaaaaaaaaaaaaaaaaaaaaaaaaaaaaaaaaaaa}
$$
$$
\phantom{aaaaaaaaaa}\le \il{0}{t}\exp\biggl\{-\il{\tau}{t}A(s) d s\biggr\}\bigl[|d_2(\tau)| + |d_3(\tau)|\bigr]d\tau, \ph 0 \le \tau \le t \le m_0 T. \eqno (3.15)
$$
Since $\biggl(\exp\biggl\{2\il{0}{t}\bigl[a_0(\tau) q_0(\tau) + a_1(\tau) q_1(\tau)\bigr]d\tau\biggr)' \ge 0, \ph t \ge 0$, from (3.13)-(3.15) we obtain
$$
|J(t)| \le \mathfrak{M}  , \phh t \in [0,m_0 T], \eqno (3.16)
$$
where
$$
\mathfrak{M}\equiv \sup\limits_{0\le\tau \le t\le m_0 T}\exp\biggl\{-\il{\tau}{t}\bigl(b_0(s) + c_0(s)\bigr) d s\biggr\}\biggl[\frac{|b_2(\tau) + c_2(\tau)| + |b_3(\tau) + c_3(\tau)|}{2a_0(\tau)} +
$$
$$
\frac{|b_2(\tau) - c_2(\tau)| +
+|b_3(\tau) - c_3(\tau)|}{2a_1(\tau)}\biggr]_0 +
$$
$$
+\max\limits_{0\le t\le m_0 T}\il{0}{t}\exp\biggl\{-\il{\tau}{t}\bigl[b_0(s) = c_0(s)\bigr]d s\biggr\}\biggl[|d_2(s)| + |d_3(s)|\biggr]d\tau < + \infty.
$$
This together with (3.12) implies that
$$
|w(t)| \le \exp\biggl\{-\il{0}{t}\bigl[b_0(\tau) + c_0(\tau)\bigr]d\tau\bigg\}|w(0)| + \mathfrak{M}, \phh 0 \le t \le m_0 T.
$$
Hence
$$
|q_0(t)| \le \exp\biggl\{-\il{0}{t}\bigl[b_0(\tau) + c_0(\tau)\bigr]d\tau\bigg\}\bigl[|q_2(0)| + |q_3(0)|\bigr] + \mathfrak{M}, n=2,3, 0\le t \le m_0 T. \eqno (3.17)
$$
Consider the integral
$$
J_0\equiv \il{0}{m_0 T}\exp\biggl\{-\il{\tau}{m_0 T}\bigl[a(s) z(s) + \widetilde{b}(s) + \widetilde{c}(s)\bigr]d s\biggr\} U(\tau) d \tau.
$$
Taking into account (3.17) one can easily show that $U(\tau)| \le c_1 |q_2(0)|^2 + c_2 |q_3(0)|^2 + c_3, \ph 0\le \tau\le m_0 T,$ where $c_1, \ph c_2, \ph c_3$ are some constants. Then
$$
|J_0|\le [c_1 |q_2(0)|^2 + c_2 |q_3(0)|^2 + c_3]\il{0}{m_0 T}\exp\bigg\{-\il{\tau}{m_0 T}(b_0(s) + c_0(s)) d s\bigg\} d\tau. \eqno (3.18)
$$
By (3.12) we have
$$
w(m_0 T)  = \exp\biggl\{-\il{0}{m_0 T}A(\tau) d\tau\biggr\}w(0)J(m_0 T). \eqno (3.19)
$$
It follows from (3.3) and from the inequality $a_0(t)q_0(t) + a_1(t) q_1(t) \ge 0, \ph t \ge 0$ that
$$
\exp\biggl\{-\il{0}{m_0 T}A(\tau) d\tau\biggr\} \le \frac{1}{e}. \eqno (3.20)
$$
Therefore $\rho \equiv 1 - \exp\biggl\{-\il{0}{m_0 T}A(\tau) d\tau\biggr\} \ne 0$ and thus fom (3.19) we obtain
$$
w(0) - w(m_0 T) = \biggl[w(0) - \frac{J(m_0 T)}{\rho}\biggr]\rho. \eqno (3.21)
$$
From (3.3) and from the inequality $a_0(t)q_0(t) + a_1(t) q_1(t) \ge 0, \ph t \ge 0$ it follows also that
$$
\exp\biggl\{-\il{0}{m_0 T}\bigl[a(\tau) z(\tau) + \widetilde{b}(\tau) + \widetilde{c}(\tau)\bigr] d\tau\biggr\} \le \frac{1}{e}. \eqno (3.22)
$$
So, we have
$$
\rho_1 \equiv 1 - \exp\biggl\{-\il{0}{m_0 T}\bigl[a(\tau) z(\tau) + \widetilde{b}(\tau) + \widetilde{c}(\tau)\bigr] d\tau\biggr\} \ne 0.
$$
Then from (3.11) we obtain
$$
z(0) - z(m_0 T) = \biggl[z(0) - \frac{J_0}{\rho_1}\biggr] \rho_1. \eqno (3.23)
$$
Let $v_0(t)\equiv v_{0,0} - i v_{0,1}(t) - j v_{0,2}(t) - k v_{0,3}(t), \ph t \ge 0$ be the solution of Eq. (1.10 with $v_0(0) = 0$. By virtue of Theorem 2.1 from the conditions 1) and 2) of the theorem it follows that $v_0(t)$ exists on $[0,+\infty)$ and $v_{0,n}(t) \ge o, \ph n=0,1, \ph t \ge 0.$ Set
$$
\xi_n \equiv \sist{sign\hskip 2pt v_{0,n}(m_0 T), \ph if \ph v_{0,n}(m_0 T) \ne 0,}{1, \phh if \ph v_{0,n}(m_0 T) = 0, \phh n = \overline{0,3}.}
$$
Let $\mathbf{v}_{\lambda,\mu} \equiv v_{\lambda,0}(t) - i v_{\lambda,1}(t) - j v_{\mu,2}(t) - k v_{\mu,3}(t), \ph t \ge 0$ be another solution of Eq. (1.1) with $v_{\lambda,n}(0) = \lambda \xi_n, \ph n= 0,1, \ph v_{\mu,n}(0) = \mu\xi_n, \ph n =2,3, \ph \lambda > 0, \ph \mu > 0.$ Since $v_{0,n}(m_0 T) \ge 0, \ph m=0,1$ in virtue of Theorem 2.1 from 1) and 2) it follows that  $\mathbf{v}_{\lambda,\mu}$ exists on $[0,+\infty)$ for every $\lambda > 0, \ph \mu > 0$. Set:
$z_\lambda(t) \equiv v_{\lambda,0}(t)- iv_{\lambda,1}(t), \ph w_\mu(t)\equiv v_{\mu,2}(t) + i v_{\mu,3}(t), \ph A_\lambda(t)  \equiv b_0(t)+ c_0(t) + 2[a_0(t) v_{\lambda,0}(t) + a_1(t) v_{\lambda,1}(t)] + i [b_1(t) - c_1(t)], V_\lambda(t) \equiv -(b_2(t) + c_2(t)) v_{\lambda,0}(t) + (b_2(t) - c_2(t)) v_{\lambda,1}(t) - i [(b_3(t) + c_3(t)) v_{\lambda,0}(t) + (b_2(t) - c_2(t)) v_{\lambda,1}(t)] - d_2(t) - i d_3(t), \ph t \ge 0$. Then by (3.21) we have
$$
w_\mu(0) - w_\mu(m_0 T) = \mu\biggl[\xi_2 + i \xi_3 - \frac{J_\lambda}{\mu\rho_\lambda}\biggr]\rho_\lambda, \eqno (3.24)
$$
where $J_\lambda \equiv \il{o}{m_0 T}\exp\biggl\{-\il{\tau}{m_0 T} A_\lambda(s) d s\biggr\} V_\lambda(\tau) d \tau, \ph \rho_\lambda \equiv 1 - \exp\biggl\{- \il{0}{m_0 T} A_\lambda(s) d s\biggr\} \ne 0.$ By virtue of (3.16) we have that for all $\lambda > 0$
$$
|J_\lambda| \le \mathfrak{M} \eqno (3.25)
$$
and by (3.30) it follows that $|Arg \hskip 2pt \rho_\lambda| < \frac{\pi}{4}$. This together with (3.24) and (3.25) implies that (see pict. 1)
$$
sign \hskip 2pt (v_{\mu,n}(0) - v_{\mu,n}(m_0 T) = \xi n, \ph n =2,3 \eqno (3.26)
$$

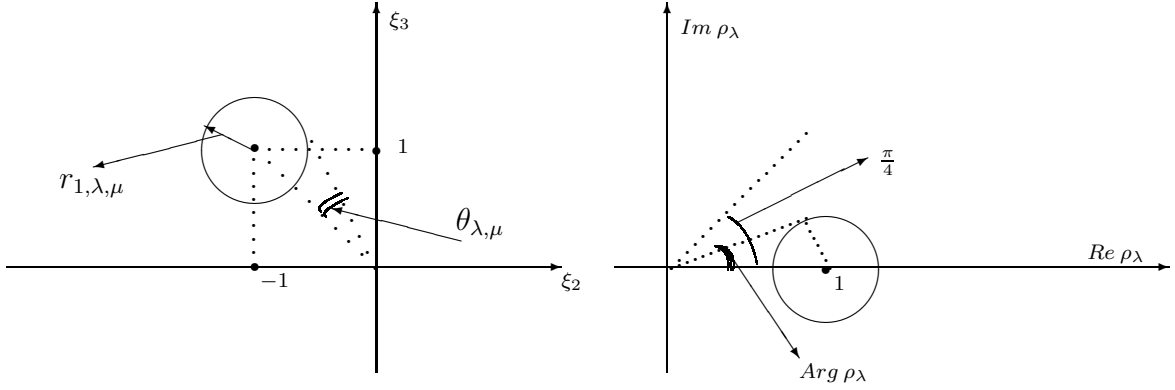
\begin{picture}(140,230)
\put(130,60){\vector(0,1){140}}
\put(-10,100){\vector(1,0){210}}
\put(240,60){\vector(0,1){140}}
\put(220,100){\vector(1,0){210}}
\put(135,193){$_{\xi_3}$}
\put(200,94){$_{\xi_2}$}
\put(138,145){$_{1}$}
\put(86,94){$_{-1}$}
\put(84,145){\circle*{3}}
\put(130,144){\circle*{3}}
\put(84,100){\circle*{3}}
\put(84,144){\circle{120}}
\put(300,99){\circle{120}}
\put(300,99){\circle*{3}}
\put(303,93){$_1$}
\put(245,190){$_{Im \hskip 2pt \rho_\lambda}$}
\put(399,105){$_{Re \hskip 2pt \rho_\lambda}$}

\multiput(128,99)(-5,5){9}{$.$}
\multiput(128,144)(-5,0){9}{$.$}
\multiput(82,100)(0,5){9}{$.$}
\multiput(128,99)(-3,6){9}{$.$}
\multiput(240,99)(3,1.1){18}{$.$}
\multiput(300,99)(-1.5,2.8){7}{$.$}
\multiput(240,99)(3,3){18}{$.$}

\put(84,144){\vector(-2,1){19}}
\put(10,130){$r_{1,\lambda,\mu}$}
\put(160,114){$\theta_{\lambda,\mu}$}
\put(290,60){$_{Arg \hskip 2pt \rho_\lambda}$}

\put(106,114){\qbezier(13,12)(2,7)(5,5)}
\put(104,116){\qbezier(13,12)(2,7)(5,5)}
\put(162,110){\vector(-4,1){49}}
\put(72,150){\vector(-4,-1){49}}

\put(256,114){\qbezier(18,-13)(15,1)(7,5)}
\put(256,112){\qbezier(9,-13)(10,-6)(4,-4)}
\put(255,112){\qbezier(9,-13)(10,-6)(4,-4)}
\put(254,112){\qbezier(9,-13)(10,-6)(4,-4)}
\put(265,103){\vector(2,-3){25}}

\put(267,117){\vector(2,1){49}}
\put(320,140){$_{\frac{\pi}{4}}$}

\put(220,-10){$Pict. 1$}

\put(40,20){$
r_{1,\lambda,\mu} \equiv \frac{J_\lambda}{\mu\rho_\lambda}, \ph |r_{1,\lambda,\mu}| \le \frac{\mathfrak{M}}{\mu|\rho_\lambda|}, \ph |r_{2, \mu}| = \frac{1}{e} < \frac{\sqrt{2}}{2}, \ph \theta_{\lambda,\mu} + Arg \hskip 2pt \rho_\lambda < \frac{\pi}{2}.
$}

\end{picture}

\vskip 50pt

By (3.23) we have
$$
z(0) - z_\lambda(m_0 T) = \lambda\biggr[\xi_0 + i\xi_1 - \frac{J_{\lambda,\mu}}{\lambda\widetilde{\rho}_\lambda}\biggr]
\widetilde{\rho}_\lambda \eqno (3.27)
$$
where $J_{\lambda,\mu} \equiv \il{0}{m_0 T}\exp\biggl\{-\il{\tau}{m_0 T}[a(s) z_\lambda(s) + \widetilde{b}(s) + \widetilde{c}(s)] d s\biggr\} U_\mu(\tau) d \tau, \ph U_\mu(t) \equiv \linebreak \widetilde{P}(t,v_{\mu,2}(t),v_{\mu,3}(t)) - i\widetilde{Q}(t,v_{\mu,2}(t),v_{\mu,3}(t)), \ph t \ge 0, \ph \widetilde{\rho}_\lambda \equiv 1 - exp\biggl\{-\il{0}{m_0 T}[a(\tau) z_\lambda(\tau) + \widetilde{b}(\tau) + \widetilde{c}(\tau)] d \tau\biggr\} \ne 0.$ By (3.3) we have
$$
|Arg \hskip 2pt \widetilde{\rho}_\lambda| <\frac{\pi}{4} \eqno (3.28)
$$
By (3.18) $J_{\lambda,\mu}$ is uniformly bounded with respect to $\lambda > 0$. From here, from (3.27) and (3.28) by analogy of (3.26) one can obtain the relations
$$
sign \hskip 2pt (v_{\lambda,n}(0) - v_{\lambda,n}(m_0 T)) = \xi_n, \ph n=0,1, \eqno (3.29)
$$
for enough large $\lambda > 0$. Assume the parameters $\lambda > 0$ and $\mu > 0$ are already chosen so large that (3.26) and (3.29) are valid. Determine the sequence of intervals $\{[\alpha_{l,n},\beta_{l,n}]_{l=0}^{+\infty}, \ph n=\overline{0,3}$ and the sequence of solutions $\{V_l(t)\}_{l=0}^{+\infty}$ of solutions of Eq. (1.1) by induction on $l$ as follows

\noindent
$
V_0(t) \equiv V_{0,0}(t) - i V_{0,1}(t) - j V_{0,2}(t) - k V_{0,3}(t), \ph t \ge 0, \ph with \ph V_{0,n}(0) = 0, \ph n=\overline{0,3}, \ph V_1(t) \equiv \linebreak V_{1,0}(t) - i v_{1,1}(t)- j V_{1,2}(t) - K V_{1,3}(t) = \mathbf{v}_{\lambda,\mu}(t), \ph t \ge 0, \ph \alpha_{0,n} = 0, \ph \beta_{0,n} = V_{1,n}(0), \ph n=0,1, \ph \alpha_{0,n}=\min\{0,V_{1,n}(0)\}, \ph \beta_{0,n} = \max\{0,,V_{1,n}(0)\}, \ph n=2,3.
$
Let for some $l$ the intervals $[\alpha_{l,n},\beta_{l,n}], \ph n= \overline{0,3}$ and the solutions $V_{l+1}(t)$ are already determined. Determine $V_{l+2}(t) \equiv V_{l+2,0}(t) - i V_{l+2,1}(t) - j V_{l+2,2}(t) - k V_{l+2,3}(t), \ph t \ge 0$ and $[\alpha_{l,n+1},\beta_{l,n+1}], \ph n= \overline{0,3}$ as follows $V_{l+2,n}(0) = \frac{\alpha_{l,n} + \beta_{l,n}}{2}\stackrel{def}{=} \gamma_{l,n}, \ph n = \overline{0,3},$
$$
\alpha_{l+1,n} = \sist{\alpha_{l,n}, \ph if \ph V_{l+2,n}(0) \ge V_{l+2,n}(m_0 T),}{\gamma_{l,n}, \ph if \ph V_{l+2,n}(0) < V_{l+2,n}(m_0 T),} \phantom{aaaaaaaaaaaaaaaaaaaaaaaaaaaaaaaaaaaaa}
$$
$$
\phantom{aaaaaaaaaaaaaaaaaaaaaaaaaaa}\beta_{l+1,n} = \sist{\beta_{l,n}, \ph if \ph V_{l+2,n}(0) < V_{l+2,n}(m_0 T,)}{\gamma_{l,n}, \ph if \ph V_{l+2,n}(0) \ge V_{l+2,n}(m_0 T),} \ph n = \overline{0,3}.
$$
Obviously $\alpha_{l,n} \ge 0, \ph n=0,1, \ph l=0,1,\dots.$ Then according to Theorem 2.1 all solutions $V_l(t), \ph l = 0,1,\dots$ exist on $[0,m_0 T]$ and, hence, all intervals $[\alpha_{l,n}, \beta_{l,n}], \ph n = \overline{0,3}, \ph l=0,1,\dots$ are determined correctly. It is also obvious that $[\alpha_{0,n},\beta_{0,n}] \supset [\alpha_{1,n},\beta_{1,n}] \supset\dots \supset [\alpha_{l,n},\beta_{l,n}] \supset \dots$ and
$$
\beta_{l,n} - \alpha_{l,n} = \frac{\beta_{0,n} -\alpha_{0,n}}{2^l}, \ph n = \overline{0,3}, \ph l=1,2,\dots. \eqno (3.30)
$$
Therefore
$$
\cap_{l=0}^{+\infty} \stackrel{def}{=} \{\gamma_n\} \not = \emptyset, \phh n = \overline{0,3}. \eqno (3.31)
$$
Show that the solution $u(t) \equiv u_0(t) - i u_1(t) - j u_2(t) - k u_3(t)$ of the system (1,1) with $u_n(0) = \gamma_n, \ph n=\overline{0,n}$ is $m_0 T$-periodic. Since $\gamma_n \ge 0, \ph n=0,1$ by virtue of Theorem 2.1 $u(t)$ exists on $[0,m_0 T]$. We must show that
$$
u(0) = u(m_0 T). \eqno (3.32)
$$
Suppose for some $n \in \{0, \ph 1, \ph 2, \ph 3\}$
$$
u_n(0) < u_n(m_0 T). \eqno (3.33)
$$
It follows from (3.30) that for every $\varepsilon > 0$ there exists $l=l(\varepsilon)$ such that
$$
0\le \gamma_n - \alpha_{l,n} \le \varepsilon, \phh 0 \le \beta_{l,n} - \gamma_n \le \varepsilon. \eqno (3.34)
$$
Let $l_0$ be chosen so that
$$
V_{l_0,n}(0) = \beta_{l,n}, \phh n = \overline{0,3}.
$$
Then from (3.25) (or (3.29)) and from the way of determination of solutions $V_l(t), \linebreak l=0,1,\dots$ is seen that
$$
V_{l_0,n}(0) \ge V_{l_0,n}(m_0 T). \eqno (3.35)
$$
From (3.34) it follows
$$
|u_n(0) - V_{l_0,n}(0)| = |\gamma_n - \beta_{l,n}| \le \varepsilon. \eqno (3.36)
$$
Since the solutions of Eq. (1.1) continuously depend on their initial values, we can chose $l(\varepsilon)$ so large that ($l_0$ depends on $l(\varepsilon)$)
$
|u_n(m_0 T) - V_{l_0,n}(m_0 T)| \le \varepsilon.
$
This together with (3.35) and (3.36) implies that $u_n(m_0 T) - u_n(0) = [u_n(m_0 T) - V_{l_0,n}(m_0 T)] + [V_{l_0,n}(m_0 T) - V_{l_0,n}(0)] + [ V_{l_0,n}(0) - u_n(0)] \le 2 \varepsilon.$ Therefore $u_n(0) \ge u_n(m_0 T)$, which contradicts (3.3). By analogy one can show that the relation $u_n(0) > u_n(m_0 T)$ leads to a contradiction. Hence $u_n(0) = u_n(m_0 T)$. Since $n$ is chosen arbitrarily from the  set $n \in \{0, \ph 1, \ph 2, \ph 3\}$ the equality (3.32) is valid.
If $D_0(t)D_1(t)\not \equiv 0$, then there exists an interval $(\alpha,\beta) \subset [0,T]$ such that
$$
D_n(t) < 0, \ph t \in (\alpha,\beta), \ph n=0,1. \eqno (3.37)
$$
Let $\widetilde{q}(t)\equiv \widetilde{q}_0(t) - i \widetilde{q}_1(t) - j \widetilde{q}_2(t) - k \widetilde{q}_3(t)$ be a solution of Eq. (1.1) with $\widetilde{q}_n(0)\ge~ 0, \linebreak n=~0,1 \ph (\widetilde{q}_n(0) \in \mathbb{R}, \ph n=2,3$). By Theorem 2.1 $\widetilde{q}(t)$ exists on $[0,+\infty)$ and
$$
\widetilde{q}_n(t) \ge 0, \ph t \ge 0, \ph n=0,1. \eqno (3.38)
$$
Show that $\widetilde{q}(t)$ is $0$-normal. By the first and the second equalities of (3.1) we can interpret $\widetilde{q}_0(t)$ and $\widetilde{q}_1(t)$ as solutions of the Riccati equations
$$
x' + a_0(t) x^2  + \{b_0(t) + c_0(t) + 2a_1(t)\widetilde{q}_1(t)\} x + P_1(t) = 0, \ph t \ge 0,
$$
$$
x' + a_0(t) x^2  + \{b_0(t) + c_0(t) + 2a_0(t)\widetilde{q}_0(t)\} x + Q_1(t) = 0, \ph t \ge 0
$$
respectively, where
$$
P_1(t) \equiv - a_0(t) \widetilde{q}_1^2(t) + (b_1(t) + c_1(t)) \widetilde{q}_1(t) - \widetilde{P}(t,\widetilde{q}_2(t),\widetilde{q}_3(t)),
$$
$$
Q_1(t) \equiv - a_1(t) \widetilde{q}_0^2(t) - (b_1(t) + c_1(t)) \widetilde{q}_0(t) - \widetilde{Q}(t,\widetilde{q}_2(t),\widetilde{q}_3(t)), \ph t \ge 0.
$$
It follows from (3.37) that $P_1(t) < 0, \ph Q_1(t) < 0, \ph t \in (\alpha,\beta)$. This together with (3.38) implies that
$$
\widetilde{q}_n(t) > 0, \ph t \in (\alpha,\beta), \ph n=0,1.  \eqno (3.39)
$$
Indeed, if $\widetilde{q}_0(t_1) = 0 \ph (\widetilde{q}_1(t_1) = 0)$ for some $t_1\in (\alpha,\beta)$ then from the relations $P_1(t_1) <~ 0, \linebreak Q_1(t_1) < 0$ it follows that $\widetilde{q}_0(t_2) < 0 \ph (\widetilde{q}_1(t_2) < 0)$ for some $t_2< t_1$, enough close to $t_1$, which contradicts (3.38). By Theorem 2.1 from (3.39) it follows that $\widetilde{q}(t)$ (in particular $q(t)$0 is $\frac{\alpha + \beta}{2}$-normal. Then by virtue of Theorem 2.2. the function $\mathbb{I}_{q, \widetilde{q}}(t) \equiv \linebreak \il{\frac{\alpha + \beta}{2}}{t} Re \hskip 2pt [a(\tau)(q(\tau) - \widetilde{q}(\tau))]d \tau, \ph t \ge \frac{\alpha + \beta}{2}$ is bounded on $[\frac{\alpha + \beta}{2}, +\infty)$. Therefore the function $I_{q,\widetilde{q}}(t), \ph t \ge 0$ is bounded on $0,+\infty)$, and by Theorem 2.2 $\widetilde{q}(t)$ (in particular $q(t)$) is $0$-normal.
The theorem is proved.

{\bf Remark 3.1}. {\it From the proof of Theorem 3.1 is seen that we can put $m_0 =1$ in it if we replace the condition 5) by the following one
$$
\il{0}{T}[b_0(\tau) + c_0(\tau)] d \tau > \ln \sqrt{2} \simeq 0.35.
$$
}
Hereafter we will assume that the functions $a(t), \ph b(t), \ph c(t)$ and $d(t)$ are continued on whole axis $\mathbb{R}$ by periodicity.
In Eq. (1.1) substitute $q \rightarrow - q, \phh t \rightarrow - t$. We come to the equation
$$
q' + q a(-t) q - b(-t) q - q c(-t) + d(-t) = 0,
$$
which together with Theorem 3.1 implies

{\bf Corollary 3.1.} {\it Let the conditions 1) - 4) of Theorem 3.1 and the condition
$$
 \il{0}{T}[b_0(\tau) + c_0(\tau)] d\tau \le 0
$$
be satisfied. Then Eq. (1.1) has a $m_0 T$-periodic solution $q^*(t) \equiv q_0^*(t) - i q_1^*(t) - j q_2^*(t) - k q_3^*(t)$ for some natural $m_0$ such that $q_n^*(t) \le 0, \ph n= 0,1, \ph t \in \mathbb{R}.$
}

Combining this result with Theorem 3.1 we obtain the following assertion

{\bf Corollary 3.2.} {\it Let the conditions 1) - 4) of Theorem 3.1 and the condition
$$
 \il{0}{T}[b_0(\tau) + c_0(\tau)] d\tau = 0 \eqno (3.40)
$$
be satisfied. Then Eq. (1.1) has  $m_0 T$-periodic solutions $q(t)\equiv q_0(t) - i q_1(t) - j q_2(t) - k q_3(t)$
and   $q^*(t) \equiv q_0^*(t) - i q_1^*(t) - j q_2^*(t) - k q_3^*(t)$ for some natural $m_0$ such that \linebreak $q_n(t) \ge 0, \ph q_n^*(t) \le 0, \ph n= 0,1, \ph t \in \mathbb{R}.$ Moreover if in addition $D_0(t)D_1(t)\not \equiv 0, \linebreak a_0(t)>~ 0$ or $a_1(t) > 0, \ph t \ge 0,$ then $q(t)$ is $0$-normal, $q^*(t)$ is $0$-extremal and the relations
$$
\limsup\limits_{t\to +\infty}\il{0}{t}Re \hskip 2pt[a(\tau)(q^*(\tau) - q(\tau))]d \tau < +\infty, \eqno (3,41)
$$
$$
\liminf\limits_{t\to +\infty}\il{0}{t}Re \hskip 2pt[a(\tau)(q^*(\tau) - q(\tau))]d \tau = - \infty \eqno (3.42)
$$
are valid.
}

Proof. Existence of $q(t)$ and $q^*(t)$, the $0$-normality and the relations $q_n(t) \ge 0, \ph q^*_n(t)\le 0, \ph t \in \mathbb{R}, \ph n=0,1$ follow immediately from Theorem 3.1 and Corollary 3.1. The $0$-extremality of $q^*(t)$ follows immediately from  Theorem 2.2 and (3.41). The relation (3.41) follows from $0$-extremality of $q^*(t)$ and from Theorem 2.2. Therefore to complete the proof of the corollary it remains to prove (3.42). It was shown in the process of proving of Theorem 3.1 that if $D_0(t)D_1(t)\not\equiv 0$, then $q_n(t) > 0, \ph t \in (\alpha,\beta), \ph n = 0,1$ for some $(\alpha,\beta)\subset[0,T]$. From here from the condition $a_0(t)>~ 0$ or $a_1(t) > 0, \ph t \ge 0,$ from the periodicity of $q(t)$ and from the inequalities $q^*_n(t) \le 0, \ph t \ge 0, \ph  n=0,1$ it follows (3.42). The corollary is proved.

{\bf Remark 3.2.} {\it The condition (3.40) is in   contrast  with the condition
$$
iii) \phantom{aaaaaaaaaa} Re \hskip 2pt [b(t) + c(t)] \le 0, \phh Im \hskip 2pt [b(t) + c(t)] \equiv 0,  \phh t \in \mathbb{R} \phantom{aaaaaaaaaaaaaaaaa}
$$
of Theorem 1.1.
}

Let us discuss now the question how we can extend the class of Riccati equations to which ca be used (indirectly but after simple  transformations) Theorem 3.1. Above we studied the case when

\noindent
I. \ph $a_0(t) \ge 0, \ph a_1(t) \ge 0, \ph a_2(t) = a_3(t) \equiv 0, \ph t\in\mathbb{R}.$

\noindent
It is not difficult to verify that the cases

\noindent
II. \ph $a_0(t) \ge 0, \ph a_1(t) \le 0, \ph a_2(t) = a_3(t) \equiv 0, \ph t\in\mathbb{R}.$

\noindent
III. \ph $a_0(t) \le 0, \ph a_1(t) \ge 0, \ph a_2(t) = a_3(t) \equiv 0, \ph t\in\mathbb{R}.$

\noindent
IV. \ph $a_0(t) \le 0, \ph a_1(t) \le 0, \ph a_2(t) = a_3(t) \equiv 0, \ph t\in\mathbb{R}.$

\noindent
can be reduced to the case I by the following simple transformations
$$
q\rightarrow \overline{q}, \phh q\rightarrow  -\overline{q}, \phh q \rightarrow - 2 \eqno (3.43)
$$
respectively. The case

\noindent
V. \ph $a_0(t) = a_1(t)\equiv 0, \ph a_2(t) \ge 0, \ph a_3(t) \ge 0, \ph t\in \mathbb{R}$

\noindent
is reducible to the case III by the transformation
$$
q \rightarrow j q. \eqno (3.44)
$$
Therefore the cases

\noindent
VI. \ph  $a_0(t) = a_1(t)\equiv 0, \ph a_2(t) \ge 0, \ph a_3(t) \le 0, \ph t\in \mathbb{R}$,

\noindent
VII. \ph  $a_0(t) = a_1(t)\equiv 0, \ph a_2(t) \le 0, \ph a_3(t) \ge 0, \ph t\in \mathbb{R}$,

\noindent
VIII. \ph  $a_0(t) = a_1(t)\equiv 0, \ph a_2(t) \le 0, \ph a_3(t) \le 0, \ph t\in \mathbb{R}$

\noindent
can be reduced to the case III by superpositions of transformations (3.43), (3,44). It is not difficult to verify that Theorem 3.1 remains valid if we replace its condition 1) by one of the following conditions (after replacement the obtained assertion can be proved by analogy of the proof of Theorem 3.1)

\noindent
IX. \ph $a_0(t) \ge 0, \ph a_2(t)\ge 0, \ph a_1(t) = a_3(t) \equiv 0, \ph t \in \mathbb{R},$

\noindent
X. \ph $a_0(t) \ge 0, \ph a_3(t)\ge 0, \ph a_1(t) = a_2(t) \equiv 0, \ph t \in \mathbb{R}.$

\noindent
To these cases can be reduced by combinations of transformations (3.43) and $q\rightarrow i q$ the following ones

\noindent
XI. \ph $a_0(t) \ge 0, \ph a_2(t)\le 0, \ph a_1(t) = a_3(t) \equiv 0, \ph t \in \mathbb{R},$

\noindent
XII. \ph $a_0(t) \le 0, \ph a_2(t)\ge 0, \ph a_1(t) = a_3(t) \equiv 0, \ph t \in \mathbb{R},$

\noindent
XIII. \ph $a_0(t) \le 0, \ph a_2(t)\le 0, \ph a_1(t) = a_3(t) \equiv 0, \ph t \in \mathbb{R},$

\noindent
XIV. \ph $a_0(t) \ge 0, \ph a_3(t)\le 0, \ph a_1(t) = a_2(t) \equiv 0, \ph t \in \mathbb{R}.$

\noindent
XV. \ph $a_0(t) \le 0, \ph a_3(t)\ge 0, \ph a_1(t) = a_2(t) \equiv 0, \ph t \in \mathbb{R}$,

\noindent
XVI. \ph $a_0(t) \le 0, \ph a_3(t)\le 0, \ph a_1(t) = a_2(t) \equiv 0, \ph t \in \mathbb{R},$

\noindent
XVII. \ph $a_1(t) \ge 0, \ph a_3(t)\ge 0, \ph a_0(t) = a_2(t) \equiv 0, \ph t \in \mathbb{R},$

\noindent
XVIII. \ph $a_1(t) \ge 0, \ph a_3(t)\le  0, \ph a_0(t) = a_2(t) \equiv 0, \ph t \in \mathbb{R},$

\noindent
IXX. \ph $a_1(t) \le 0, \ph a_3(t)\ge 0, \ph a_0(t) = a_2(t) \equiv 0, \ph t \in \mathbb{R},$

\noindent
XX. \ph $a_1(t) \le 0, \ph a_3(t)\le 0, \ph a_0(t) = a_2(t) \equiv 0, \ph t \in \mathbb{R},$

\noindent
XXI. \ph $a_1(t) \ge 0, \ph a_2(t)\ge 0, \ph a_0(t) = a_3(t) \equiv 0, \ph t \in \mathbb{R},$

\noindent
XXII. \ph $a_1(t) \ge 0, \ph a_2(t)\le 0, \ph a_0(t) = a_3(t) \equiv 0, \ph t \in \mathbb{R},$

\noindent
XXII. \ph $a_1(t) \le 0, \ph a_2(t)\ge 0, \ph a_0(t) = a_3(t) \equiv 0, \ph t \in \mathbb{R},$

\noindent
XXIV. \ph $a_1(t) \le 0, \ph a_2(t)\le 0, \ph a_0(t) = a_3(t) \equiv 0, \ph t \in \mathbb{R}.$

\noindent
We see that by simple transformations the area of application of Theorem 3.1 is enough wide. The following  approach shows that this are (the cases I-XXIV) can be radically extended. Let $\lambda(t)$ be a quaternionic-valued continuously differentiable function on $\mathbb{R}$ such that $\lambda(t)\ne 0, \ph t \in \mathbb{R}$. Consider the Riccati equation
$$
q' + q\lambda(t)a(t)\lambda(t)q + b(t)q + q c(t) + d(t) = 0, \phh t\in \mathbb{R},
$$
where $a(t)$ is the same as in Theorem 3.1.
Multiply both sides (at left and at right) of this equation by $\lambda(t)$. Taking into account the equality
$$
\lambda(t) q' \lambda(t) = (\lambda(t) q \lambda(t))' - \lambda'(t) q \lambda(t) - \lambda(t) q \lambda'(t), \ph t \in \mathbb{R}
$$
we obtain
$$
v' + v a(t) v + (b(t) - \lambda'(t))\lambda^{-1}(t) v +  v\lambda^{-1}(t)(c(t) - \lambda'(t))  + \lambda(t) c(t) \lambda(t), \phh t\in \mathbb{R},
$$
where $v\equiv \lambda(t) q \lambda(t), \ph  t\in \mathbb{R}.$

{\bf Remark 3.3.} {\it Unlike of Eq. (1.1) for which the components of the function $a(t)$ do not change sign, the components of the function $\lambda(t) a(t) \lambda (t)$ can change signs.}

\vskip 120pt

\centerline{\bf References}

\vskip 20pt

\noindent
1. P. Wilzinski,  Quaternionic-valued differential equations. The Riccati equations. \linebreak\phantom{aa}  Journal of Differential Equations, vol. 247. pp. 2167 - 2187, 2009.

\noindent
2. J. Campos, J. Mavhin. Periodic solutions of quaternionic-valued ordinary differential \linebreak\phantom{aa} equations. Annali di Mathematica, vol. 185, pp. 109 - 127, 2006.

\noindent
3. J. D. Gibbon, D. D. Holm, R. M. Kerr and I. Roulstone. Quaternions and periodic \linebreak\phantom{aa}  dynamics in the Euler fluid equations. Nonlinearity, vol. 19, pp. 1962 - 1983, 2006.

\noindent
4. H. Zoladek, Classification of diffeomorphisms of $\mathbb{S}^4$ induced by quaternionic Riccati \linebreak\phantom{aa}  equations with periodic coefficients. Topological methods in Nonlinear Analysis. Journal \linebreak\phantom{aa}  of the Juliusz Shauder Center, vol. 33. pp. 205 - 2015, 2009.

\noindent
5. J. Mathews, The quaternionic Structure of the Equations of Geophisical Fluid Dinamics. \linebreak\phantom{aa}  A thesis. submitted for the Degree of Doctor of Philosophy. School of Mathematics, \linebreak\phantom{aa} Metrology and Physics, 2006, 197 pages.

\noindent
6. G. A. Grigorian, Global solvability criteria for quaternionic Riccati equations. Archivum \linebreak\phantom{aa}   Mathematicum, Tomus 57 (2021), pp. 83–99.

\noindent
7. G. A. Grigorian, Properties of solutions of quaternionic Riccati equations.  Archivum \linebreak\phantom{aa}   Mathematicum. In pront.

\end{document}